\documentclass{gtart}
\usepackage{amsfonts} 
\usepackage{epsf}
\let\Bbb\mathbb

\input gtmonout
\volumenumber{2}
\volumeyear{1999}
\volumename{Proceedings of the Kirbyfest}
\pagenumbers{23}{34}
\papernumber{2}
\received{27 January 1999}\revised{15 June 1999}
\published{17 November 1999}

\def\bd{{\partial}}
\def\R{{\Bbb R}}
\def\Z{{\Bbb Z}}
\nocolon

\newtheorem{theorem}{Theorem}
\newtheorem{lemma}[theorem]{Lemma}
\newtheorem{conjecture}{Conjecture}

\begin{document}

\title{Cosmetic Surgery on Knots}
\author{Steven A Bleiler\\Craig D Hodgson\\Jeffrey R Weeks}
\shortauthors{Steven A Bleiler, Craig D Hodgson and Jeffrey R Weeks}
\begin{abstract}
This paper concerns the Dehn surgery
construction, especially
those Dehn surgeries leaving the manifold unchanged.
In particular, we describe an oriented 1--cusped hyperbolic 3--manifold
$X$ with a pair of slopes $r_1$, $r_2$ such that the Dehn filled
manifolds $X(r_1)$, $X(r_2)$ are oppositely oriented copies of the
lens space $L(49,18)$, and there is no homeomorphism $h$ of $X$ such that
$h(r_1) = h(r_2)$.
\end{abstract}
\asciiabstract{This paper concerns the Dehn surgery
construction, especially
those Dehn surgeries leaving the manifold unchanged.
In particular, we describe an oriented 1-cusped hyperbolic 3-manifold
X with a pair of slopes r_1, r_2 such that the Dehn filled
manifolds X(r_1), X(r_2) are oppositely oriented copies of the
lens space L(49,18), and there is no homeomorphism h of X such that
h(r_1)=h(r_2).}

\keywords{Dehn surgery, Dehn filling, hyperbolic knots}
\primaryclass{57N10}\secondaryclass{57M25, 57M50}
\maketitle

\section{Introduction}

This paper concerns the Dehn surgery construction \cite{D}, in
particular, those Dehn surgeries leaving the manifold unchanged,
a phenomenon we call {\it cosmetic} surgery.  In what
follows it will be useful to consider the Dehn surgery
construction in terms of the slightly more general notion of
Dehn filling, recalled below.

Let $X$ be an oriented 3--manifold with torus boundary $\partial X = T^2$, and
let $r$ be a slope on $\partial X$, that is, an isotopy class of
unoriented, simple, closed
curves in $T^2$.  Define $X(r)$, the {\it r--Dehn filling} on $X$,
as $ X \cup
(B^2 \times S^1)$
where the boundaries are identified by a homeomorphism taking
$\partial B^2 $ to the slope $r$.
The connection with Dehn surgery is regained by considering
$X(r_2)$ as a surgery on $M = X(r_1)$.  Two Dehn fillings
are considered {\it cosmetic} if there is a
homeomorphism $h$ between $X(r_1)$ and $X(r_2)$.
The Dehn fillings are {\it truly cosmetic} if the
homeomorphism $h$ is orientation preserving,
and {\it reflectively cosmetic} if $h$ is orientation reversing.
A pair of fillings $X(r_1)$ and $X(r_2)$ may
be both truly cosmetic and reflectively cosmetic.

One's intuition might say that cosmetic surgeries must be rare, but examples
are easy to find.  As described in Rolfsen's beautiful book \cite[Chapter 9]{R}
 reflecting in the plane of the paper shows that $p/q$ and
$-p/q$ surgeries on an amphicheiral knot in the 3--sphere yield oppositely
oriented copies of the same manifold, and performing $1/n$ surgery on
the unknot always produces a consistently oriented 3--sphere.  See Figure 1.

\begin{figure}[ht!]
\centerline{\epsfxsize.9\hsize\epsfbox{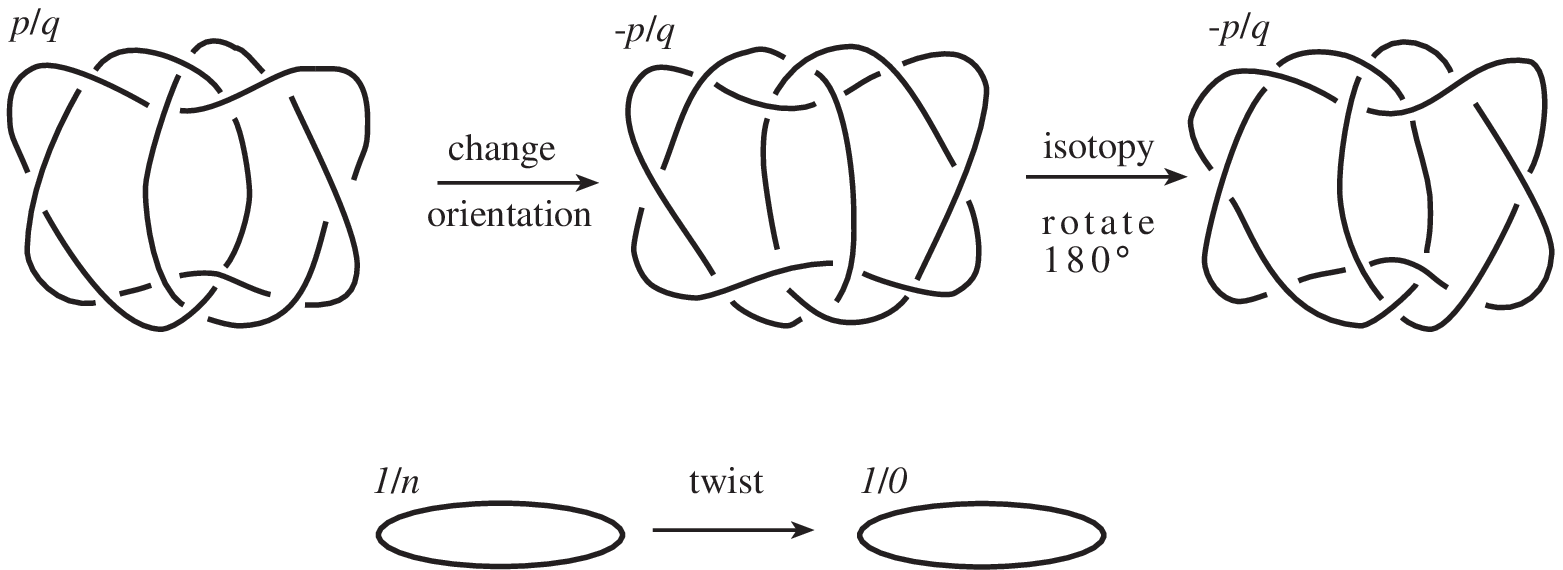}}
\caption{}
\end{figure}

In both these cases, however, there exists a homeomorphism of
the torally bounded knot exterior $X$ which takes one surgery
slope to the other.  We call two slopes {\it equivalent} if
such a homeomorphism exists.
Truly cosmetic (resp.\  reflectively cosmetic)
Dehn fillings $X(r_1)$ and $X(r_2)$
are considered {\it mundane} when there is an orientation
preserving (resp.\ orientation reversing) homeomorphism of $X$
taking $r_1$ to $r_2$.
Cosmetic Dehn fillings which are not mundane are {\it exotic}.

There is one easy way to find exotic cosmetic fillings:
swapping the sides of a Heegaard splitting of certain  lens
spaces can produce
true and reflective exotic cosmetic
surgeries on the unknot.
In particular, $p/q$ and $p/q'$
surgeries yielding the
lens spaces
$L(p,q)$ and $L(p,q')$ where
$q q' \equiv 1 \bmod p$ and $q' \ne q \bmod p$
perform this happy trick,
for example $17/2$ and $17/9$ surgeries.
So we require the space $X$ not to be homeomorphic
to $B^2 \times S^1$ in our definition of {\it exotic} cosmetic fillings.

Now things are much harder.  Indeed, recent results in surgery
theory suggest that examples of exotic cosmetic surgeries
are few and far between.  For example, the
solutions to the knot complement problems in $S^3$ \cite{GL} and
in $S^2 \times S^1$ \cite{G1},
phrased in the language here,
state that the 3--sphere $S^3$ and the manifold $ S^2 \times S^1$ never
arise via
cosmetic surgery.  So one thinks of cosmetic surgeries as occurring on
knots in more general manifolds.

Similar results hold when the first Betti number of the
manifold is positive.  For example, Boileau, Domergue, and Mathieu \cite{BDM}
showed in 1995 that if $X$ is irreducible and the core of the surgered
solid torus is
homotopically trivial in such an $X(r_1)$, then for $r_2$ distinct from
$r_1$, the
manifold $X(r_2)$ is never even simple homotopy equivalent to $X(r_1)$.  An
extension and sharpening of this result was given by  M Lackenby \cite{L}
shortly thereafter.  Assume that $X$ is irreducible and atoroidal and, as
before, that the
core of the surgered torus is homotopically trivial in $X(r_1)$ with first
Betti number positive. Then
Lackenby's theorem shows that if at least one of the slopes $r_2$ and $r_3$ has
a sufficiently high geometric intersection number with $r_1$, $X(r_2)$ and
$X(r_3)$ are orientation preserving homeomorphic if and only if $r_2 =
r_3$, and $X(r_2)$ and $X(r_3)$ are orientation reversing homeomorphic if
and only if the surgery core is amphicheiral  and $r_2 = - r_3$.

The assumption of homotopic triviality is important
here, as illustrated by
certain Seifert fibre spaces.
As first shown by Mathieu \cite{M}, $9/1$ and
$9/2$ surgeries on the right hand trefoil yield oppositely oriented copies
of the same Seifert fibre space. The slopes are inequivalent
since they have
different distances from the meridian, and any homeomorphism of the
knot exterior must take the meridian to itself.
The existence of Seifert fibred cosmetic
surgeries has a particularly nice picture in the Kirby calculus
\cite[Chapter 9]{R} and depends on the existence of an exceptional
fibre of index 2.
The key property of such an exceptional fibre is that after orientation
reversal, the type of an exceptional fibre of index 2 can be restored via a
twist.
See Figure 2 which begins with a surgery description of the
trefoil knot $K$ in $S^3$. (This comes from the Seifert fibration of $S^3$
over $S^2$ with a trefoil knot as regular fibre and two exceptional fibres of
order 2 and 3
 --- see Seifert \cite[sections 3,11]{seif}, Montesinos \cite[chapter 4]{Mo1}.)

\begin{figure}[ht!]
\centerline{\epsfxsize.9\hsize\epsfbox{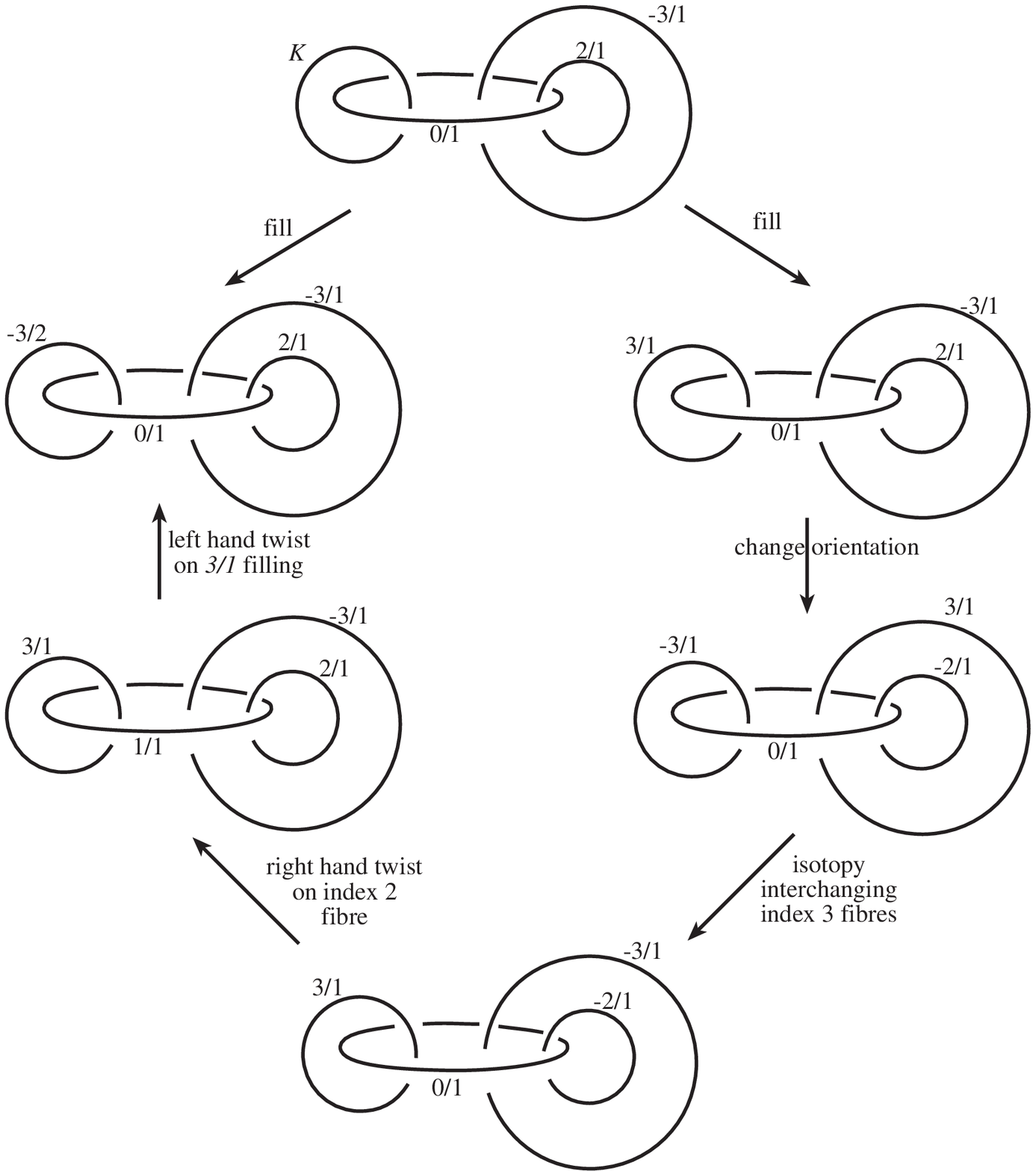}}
\caption{}
\end{figure}

The same construction gives examples of exotic cosmetic surgeries on Seifert
fibre spaces with positive first Betti number,
starting with surgery descriptions of Seifert fibred spaces over
higher genus surfaces as in \cite[figure 12, p. 146]{Mo1}.

This observation
was also used by Rong \cite{Ro} to classify the cosmetic surgeries where
$X$ carries a Seifert fibration and the surgery yields a Seifert fibre
space.   One consequence of the classification is that all the
exotic cosmetic surgeries on Seifert fibre spaces are reflective.
Another consequence of this technique is the fact that if a Seifert fibred
$X$ admits a pair of cosmetic surgeries yielding a Seifert fibre space, it
admits an infinity of such pairs.

\section{Cosmetic surgery on hyperbolic manifolds}

 From Seifert fibred cosmetic surgeries it is fairly easy to construct
examples of cosmetic surgeries on graph manifolds, so the question arises
as to whether there are any hyperbolic examples.  Such examples are
unexpected because of the general theory of 1--cusped hyperbolic manifolds.
For example, a theorem similar to M Lackenby's but for homotopically
non-trivial knots and showing there can be only finitely many cosmetic
surgeries on a hyperbolic knot is given below. Throughout this paper
hyperbolic 3--manifolds are assumed to be {\it complete} and of
{\it finite volume}.

\begin{theorem} Let $X$ be a $1$--cusped, orientable hyperbolic 3--manifold.
Then there exists a finite set of slopes $E$ on $\bd X$, such
that if $r_1$ and $r_2$ are distinct slopes outside $E$, $X(r_1)$ and
$X(r_2)$ homeomorphic implies that there exists an orientation reversing
isometry
$h$ of $X$ such that $h(r_1) =  r_2$.  In particular, if $X$ is a classical
knot complement, then the knot is amphicheiral and $r_1= - r_2$.
\end{theorem}

\noindent{\bf Proof}\qua By Thurston's theory of hyperbolic Dehn surgery
[T1], we can choose $E$ so that each filling outside
$E$ gives a hyperbolic manifold in which the
surgery core circle is isotopic to the {\sl unique} shortest closed geodesic.
Assume $r_1, r_2$ are outside $E$ and $X(r_1)$ is homeomorphic to
$X(r_2)$. Then by Mostow rigidity there is an isometry 
$X(r_1) \to X(r_2)$ taking the core geodesic $C_1$ to core geodesic $C_2$.
This restricts to a homeomorphism of $X$ taking $r_1$ to $r_2$.
The following lemma then completes the proof.\qed

\begin{lemma} Let $X$ be a
$1$--cusped, orientable 
hyperbolic 3--manifold. If $h \co  X \to X$
is a homeomorphism which changes the slope of some peripheral
curve, then $h$ is orientation reversing. If $X$ is a classical
knot complement, then $h$ takes each slope $r$ to $-r$.
\end{lemma}

\noindent{\bf Proof}\qua
By Poincar\'e duality, the map
induced by inclusion
$H_1(\bd X; \R) \to H_1(X;\R)$
has a 1--dimensional kernel $K$.  The hyperbolic structure on the
interior of $X$ gives a natural Euclidean metric on
$H_1(\bd X; \R) \cong \R^2$, defined up to similarity. By Mostow Rigidity,
$h$ is homotopic to an isometry, hence $h_* \co  H_1(\bd X; \R) \to H_1(\bd X;
\R)$ is an isometry. Further, $K$ and its orthogonal complement are preserved
by $h_*$ so each is a $+1$ or $-1$ eigenspace for $h_*$. It follows that if
$h$ is orientation preserving then $h_* = \pm {\rm~identity}$, so $h$
does not change the slope of any peripheral curve.

If $X$ is a classical knot complement, then $h$ preserves
both the longitude and meridian up to sign (using the solution
to the knot complement problem \cite{GL}).
Hence, if $h$ is orientation reversing, each slope is taken
to its negative.\qed

\bigskip

The theorem strongly suggests that finding cosmetic surgeries on
hyperbolic manifolds that yield hyperbolic manifolds is hard, perhaps
impossible.  Instead one looks to find cosmetic surgeries that yield
non-hyperbolic manifolds, such as the lens spaces.  Another reason to look
here is that there are cosmetic surgeries in the solid torus.  In fact,
Berge \cite{Be} and Gabai (\cite{G1}, \cite{G2}) classify those knots in $B^2
\times S^1$ which have inequivalent slopes which fill to $B^2 \times S^1$.

 From their classification, these knots are certain 1--bridge braids and,
with a unique exception where there are three such slopes,  these knots
have exactly two slopes which yield a solid torus when filled.  The
meridian of this new solid torus, however, is quite different from the
meridian of the original solid torus. Indeed, as originally shown by C
Gordon \cite{Go1}, if one performs $p/q$--surgery on a knot in a solid torus
and again obtains a solid torus, the meridian of this new solid torus is
given by the slope $p/( k^2 q)$, where $k$ is the winding number of the
original
knot in the solid torus.  Now, from such a knot in a solid torus it is easy to
construct a torally bounded 3--manifold having a pair of slopes yielding lens
spaces by attaching a solid torus to the outside of the solid torus in which
the knot lies.  The lens spaces produced are determined by the relation of
this last attaching curve to the meridians described above.  This has led
to a conjectural classification of the fillings on hyperbolic knot complements
in the $3$--sphere which yield lens spaces, see \cite{Go2}.

Of course, the purpose here is to construct examples where filling on the
appropriate slopes yields homeomorphic lens spaces.  This requirement
produces a severe number theoretic obstruction, namely, that the outside
attaching slope $r$ be simultaneously a $(p,q)$ and $(p,q')$ curve with
$qq' \equiv \pm 1 \bmod p$ or $q \equiv \pm q' \bmod p$
with respect to our two meridians. (See \cite{Reid} or \cite{Brody}.)
The number
$p$ is, of course, the geometric intersection number of $r$ and the
meridians, and is commonly referred to as the {\it distance} between
the two slopes.  This already puts powerful restrictions on the allowable
slopes as given a pair of slopes $r_1=a_1/b_1$ and $r_2=a_2/b_2$, there are
exactly two slopes  equidistant from $r_1$ and $r_2$.  For this one notes
that the geometric intersection number is just $|a_1 b_2 -b_1 a_2|$,
so if $r=c/d$ is
equidistant from $r_1$ and $r_2$ then either $a_1 d-b_1 c = a_2 d-b_2 c$ or
$a_1 d -b_1 c =b_2 c - a_2 d$.  In the former case it follows that
$c/d=(a_1-a_2)/(b_1-b_2)$ and in the latter that $c/d=(a_1+a_2)/(b_1+b_2)$.

\medskip

The precise lens spaces produced are then determined by the following.

\begin{lemma}
Let $M$ be the manifold obtained by
$a/b$ and $c/d$ Dehn fillings on the two boundary components
of $T^2 \times [0,1]$,
where $\gcd(a,b)=\gcd(c,d)=1$. Let $a^*,b^*$ be integers
such that $a^* b - b^* a = 1$. Then $M$ is the lens space
$L(p,q)$ where $\displaystyle{{p \over q} = {bc-ad \over a^*d-b^* c}}$.
\end{lemma}

\noindent{\bf Proof}\qua The linear automorphism of $T^2$
given by the matrix
$${\left[\matrix{ a^* & a \cr b^* & b \cr}\right]}^{-1} =
\left[\matrix{ b & -a \cr -b^* & a^* \cr}\right]$$
takes the slopes $\displaystyle{ a\over b}$ and $\displaystyle{c \over d}$ to
$\displaystyle{0\over 1}$ and $\displaystyle{bc-ad \over a^* d-b^* c}$. \qed

\medskip

A further restriction on $r$ arises from hyperbolicity.  Considering
the complements of these braids in the solid torus as link exteriors, one
notes that of the six types listed in the Berge--Gabai classification only
four are atoroidal and acylindrical, and hence by Thurston (\cite{T2},
\cite{T3}) correspond
to 2--cusped hyperbolic manifolds. If one wishes to construct hyperbolic
examples by the above construction, the slope $r$ must lie in the
hyperbolic region of the Dehn surgery space of these manifolds.  After
filling on $r$, there is also the question as to whether the slopes $r_1$
and $r_2$ are indeed inequivalent.

Both these points are exemplified by a family of braids for which the
construction gives distinct fillings yielding homeomorphic lens spaces.
These examples arise from
the Berge braids of ``Type IV" and the Pythagorean triples of the
form $(s,t,u) = (2k+1, 2k(k+1), 2k^2 + 2k +1)$.
In particular, denoting by $W_n$ the
product of the first $n-1$ standard braid generators, the braids
$W_s^{1} W_t^{-s}$
give rise to a family of torally bounded
3--manifolds each with one pair of slopes yielding
oppositely oriented copies of
$L(u + s, s + 2) = L(2 (k+1)^2, 2k+3)$
and another  pair of slopes yielding
oppositely oriented copies of $L(u - s, s - 2)= L(2k^2, 2k-1)$.
However, it turns out that each of these torally bounded
3--manifolds has the structure of an amphicheiral graph manifold.
This was initially suggested by computer calculations via SnapPea
[We], then confirmed, after using lots of string,
 by a direct analysis
via the Montesinos trick (\cite{Mo}, \cite{Bl}). So
these examples are neither hyperbolic nor surgeries on inequivalent slopes.

\section{A surprising example}

The preceding discussion makes the following example even more striking.
Our construction will produce an oriented 1--cusped hyperbolic manifold
$X$ with an exotic pair of reflective cosmetic Dehn fillings.
The filled manifolds $X(r_1)$ and $X(r_2)$ will be oppositely
oriented copies of the lens space $L(49,18)$, and the cusped
manifold $X$ will be presented naturally as a knot exterior in
$S^2 \times S^1$.

Begin with the now infamous 1--bridge braid
$W_3^{-1} W_7^3$ in a
solid torus. This braid and its mirror image are the unique 1--bridge braids
with three
distinct slopes which yield a solid torus when filled (see Figure 3).
When the
solid torus $T$ in which the braid lies is considered to be a standard
torus in the 3--sphere $S^3$, these special slopes are given by $1/0$, $18/1$,
and  $19/1$.  Filling on these slopes produces solid tori whose meridians are
represented by the
$1/0$, $18/49$, and  $19/49$ slopes on the boundary of $T$.

\begin{figure}[ht!]
\centerline{\epsfxsize.5\hsize\epsfbox{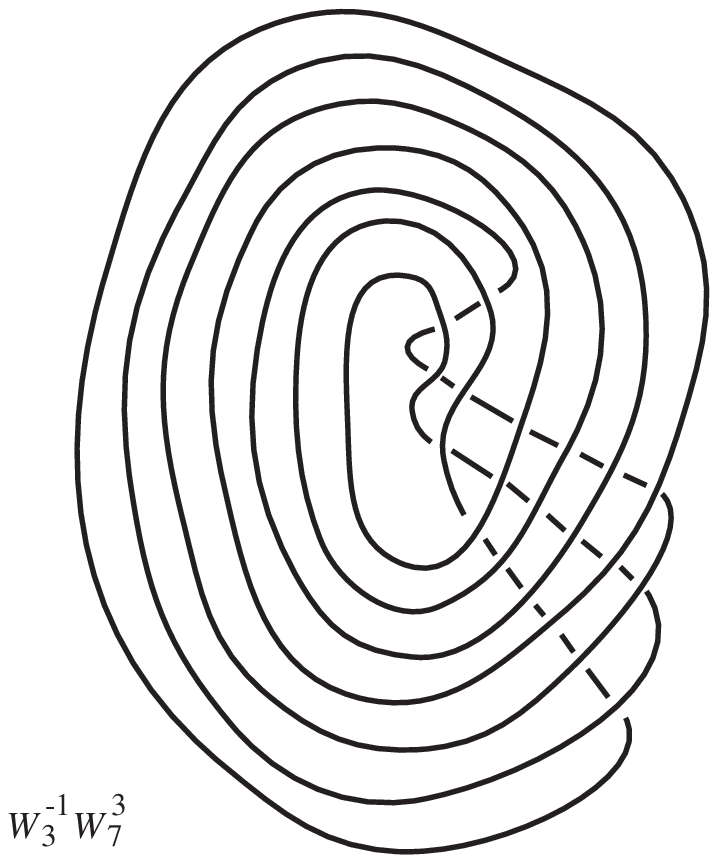}}
\caption{}
\end{figure}

For the pair $r_1=18/49$ and $r_2=19/49$ it follows from our discussion
above that the slopes $1/0$ and $37/98$ are equidistant to $r_1$ and
$r_2$.  Of course, attaching a solid torus to the $p/q$ curve on a
standard solid torus in $S^3$ is the same as performing $q/p$ surgery on
the core of the complementary solid torus.  Choosing $1/0$ then, our braid
becomes a knot in  $S^2 \times S^1$ representing $7$ in
${\pi_1} \cong \Z$.
An analysis by SnapPea \cite{We} shows that the exterior of this knot is
hyperbolic and that there is no homeomorphism
taking the slope $18/1$ to $19/1$.
(Alternatively, the Montesinos trick can be used to show this knot exterior
to be atoroidal and acylindrical.)

So if filling on these slopes does
indeed produce homeomorphic manifolds, we have the example we seek.  But
now an easy exercise, using Lemma 3 or the Kirby calculus (see Figure 4), shows
that filling on these slopes produces the manifolds $L(49,-18)$ and
$L(49,-19)$.

\begin{figure}[ht!]
\centerline{\epsfxsize.9\hsize\epsfbox{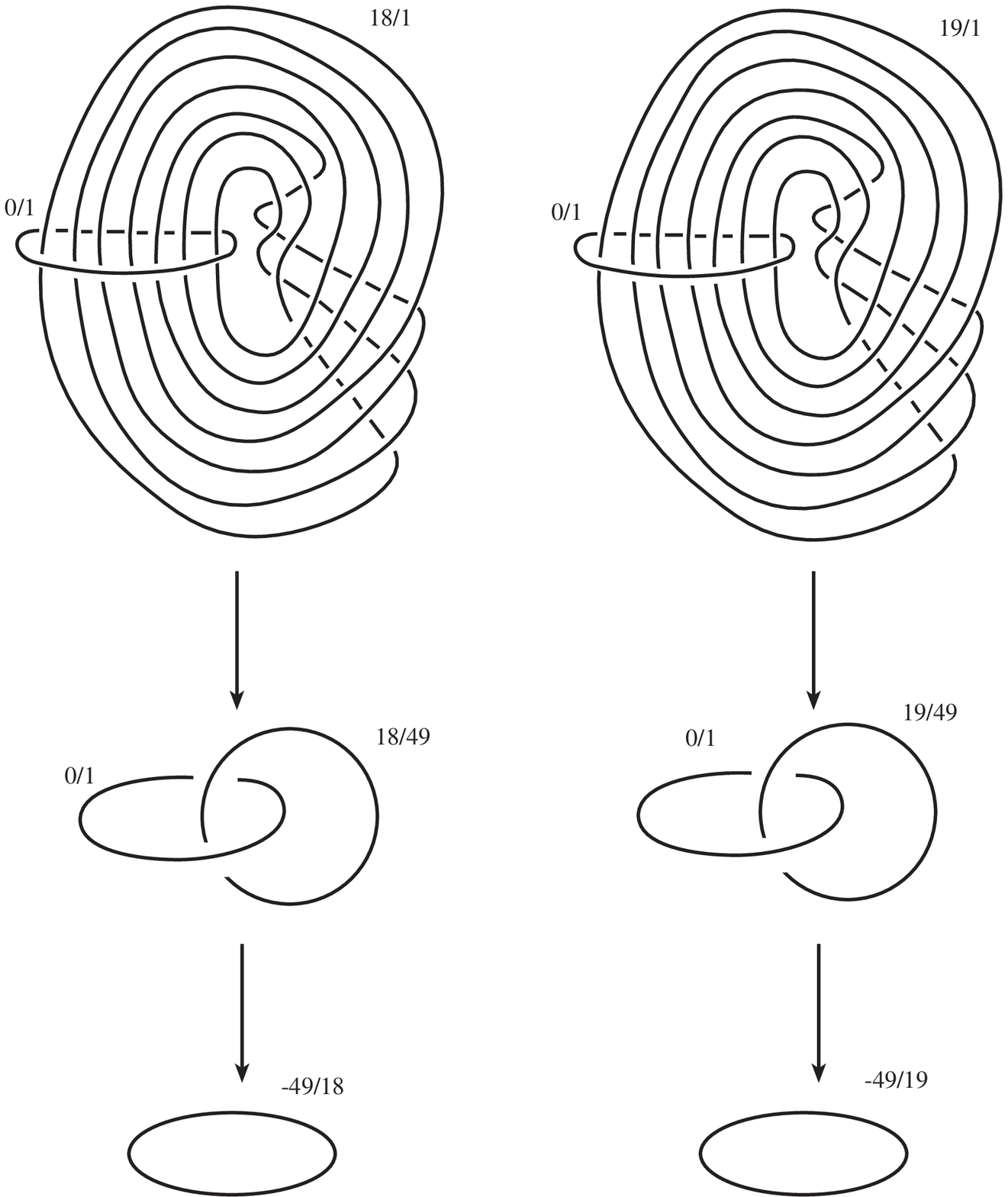}}
\caption{}
\end{figure}

Since $(-18) \cdot (-19)=342$ and $7 \cdot 49=343$,
the classification of lens spaces
shows that one has obtained oppositely oriented copies of the lens space
$L(49,18)$.

\smallskip
\noindent{\bf Remark}\qua In SnapPea's notation, these cosmetic surgeries are
given by the $(0,1)$ and $(1,1)$ fillings on the cusped hyperbolic manifold
$m172$.

\section{Concluding Remarks}

To date this is the only known example of exotic cosmetic surgery
on a hyperbolic knot exterior.  As it and all the Seifert fibred
examples are reflective, it seems appropriate to make
the following conjecture, given in problem 1.81 of Kirby's
problem list \cite{K}.

\begin{conjecture}[Cosmetic surgery conjecture] \rm
Exotic cosmetic surgeries are never truly cosmetic.
\end{conjecture}

An equivalent form of this conjecture is the

\begin{conjecture}[Oriented\, knot\, complement\, conjecture] \rm
 If $K_1$ and $K_2$ are\break
knots in a closed, oriented 3--manifold $M$ whose complements are homeomorphic
via an orientation-preserving homeomorphism, then there exists an
orientation-preserving homeomorphism of $M$ taking $K_1$ to $K_2$.
\end{conjecture}

We close with two further conjectures and a comment.  Our earlier theorem
suggests the following:

\begin{conjecture}\rm  Cusped hyperbolic manifolds admit no
cosmetic fillings, true or reflective, yielding hyperbolic manifolds.
\end{conjecture}

\begin{conjecture} \rm  Closed geodesics in a hyperbolic 3--manifold
are determined by their complements (even allowing orientation
reversing homeomorphisms).
\end{conjecture}

\smallskip
\noindent{\bf Remarks}\qua (3) $\Rightarrow$ (4) but they are not equivalent as it
may happen that the core of one of the surgeries is not isotopic to a closed
geodesic.
Some evidence for these conjectures has been provided by a computer search:
No manifold in the Hodgson--Weeks census of 11,031 low-volume
closed, orientable  hyperbolic 3--manifolds is obtained by two inequivalent
fillings on a
manifold in the Hildebrand--Weeks census of 4,815 cusped, orientable
hyperbolic 3-manifolds
triangulated by at most 7 ideal simplices. (All of these hyperbolic
manifolds are
incorporated in SnapPea \cite{We}.)

\smallskip
\noindent{\bf Question}\qua  Does there exist a pair of exotic cosmetic fillings
which are simultaneously true and reflective?

\smallskip
One can also ask the above questions with homeomorphism replaced by homotopy
equivalence or simple homotopy equivalence (see \cite{BDM} for a version which
is different because the knots are null-homotopic).
We note that there is a hyperbolic knot exterior in a lens space with a pair of
slopes which yield non-homeomorphic but orientation preserving homotopy
equivalent lens spaces.  It is obtained by replacing the $1/0$ slope in
our construction above by $37/98$.

\rk{Acknowledgements}
Steven Bleiler is partially supported by a grant from the
Portland
State University Faculty Development Program. 
Craig Hodgson is partially supported by grants from the Australian
Research Council.  Jeffrey Weeks is
partially supported by US National Science Foundation grant DMS--9803362.

\bigskip
{\small \parskip 0pt \leftskip 0pt \rightskip 0pt plus 1fil \def\\{\par}
\sl Department of Mathematics, 
Portland State University\\Portland, OR 97207-0751, USA\\
\smallskip
Department of Mathematics, University of Melbourne\\Parkville, Victoria 3052,
Australia\\
\smallskip 15 Farmer Street, Canton,  NY, USA\\
\medskip
\rm Email:\stdspace\tt steve@mth.pdx.edu\stdspace cdh@ms.unimelb.edu.au\qua
weeks@northnet.org\par}
\recd

\end{document}